\documentclass{amsart}
\usepackage{graphicx}
\usepackage{amsfonts}
\usepackage{eucal}
\usepackage{amscd}
\usepackage{amssymb}
\usepackage{xypic}
\xyoption{all}
\DeclareFontFamily{OT1}{pzc}{}
\DeclareFontShape{OT1}{pzc}{m}{it}%
             {<-> s * [1.200] pzcmi7t}{}
\DeclareMathAlphabet{\mathscr}{OT1}{pzc}%
                                 {m}{it}

\newcommand{\alg}[1]{\boldsymbol{#1}}
\newcommand{\scheme}[1]{\boldsymbol{\underline{#1}}}
\newcommand{\mc}[1]{\mathcal #1}
\newcommand{\sheaf}[1]{\mathcal{#1}}
\newcommand{\ff}[1]{\mathbf{#1}}

\newcommand{\rats}{\mathbb Q}

\newcommand{\complex}{\mathbb C}

\newcommand{\OO}{\mathcal{O}}

\newcommand{\Build}{\mathcal{B}}

\newcommand{\isom}{\cong}

\newtheorem{thm}{Theorem}[section]
\newtheorem{cor}[thm]{Corollary}

\newtheorem{prop}[thm]{Proposition}
\theoremstyle{definition}
\newtheorem{defn}[thm]{Definition}

\numberwithin{equation}{section}

\begin{document}

\title{Depth Zero Representations of Nonlinear Covers of $p$-adic Groups}%
\author{Tatiana K. Howard and Martin H. Weissman}%
\address{Tatiana K. Howard, Dept. of Mathematics, University of Michigan, Ann Arbor, MI 48109 \\ Martin H. Weissman, Dept. of Mathematics, University of California, Santa Cruz, CA 95064}%
\email{tkhoward@umich.edu, weissman@ucsc.edu}
\begin{abstract}
We generalize the methods of Moy-Prasad, in order to define and study the genuine depth zero representations of some nonlinear covers of reductive groups over $p$-adic local fields.  In particular, we construct all depth zero supercuspidal representations of the metaplectic group $Mp_{2n}$ over a $p$-adic field of odd residue characteristic.
\end{abstract}
\maketitle
\tableofcontents

\section{Introduction}

Supercuspidal representations play a crucial role in the representation theory of reductive groups over $p$-adic fields.  Beginning with results of Jacquet \cite{Jac}, every irreducible admissible representation is either supercuspidal or can be obtained (as a constituent) from parabolically inducing a supercuspidal representation.  Canonical splittings of covers of unipotent groups allow the results of Jacquet to be extended to nonlinear covers of reductive groups.  Thus supercuspidal representations are important for the representation theory of nonlinear covers of reductive groups as well.

In \cite{MP1} and \cite{MP2} Moy and Prasad introduce another invariant of representations of $p$-adic groups, called {\em depth}, which is a nonnegative rational number.   The irreducible, depth zero supercuspidal representations are most foundational, as they are used to construct almost all supercuspidal representations of positive depth (the construction of J.K.Yu \cite{Yu2} and exhaustion theorem of Ju-Lee Kim \cite{Kim}).

In this work, we define and study the irreducible genuine depth zero representations of ``tame'' nonlinear covers of reductive $p$-adic groups, where ``tameness'' refers to the assumption that the kernel of the cover has cardinality prime to $p$.  In particular, we cannot address double covers (e.g., the usual metaplectic groups) in residual characteristic two.  Under the tameness assumption, we construct all genuine depth zero supercuspidal representations of nonlinear covers $\tilde G$ of reductive $p$-adic groups, by extending then inducing irreducible cuspidal representations of related reductive groups over finite fields.  The irreducibility of such induced representations yields its supercuspidality, by an argument going back to Mautner (cf. Theorem 9.1 of \cite{Mau}).

Despite some work of Pan \cite{Pan1}, depth zero representations of nonlinear covers of reductive $p$-adic groups (e.g., the metaplectic group) have not been studied in detail.  It is the purpose of this paper to extend the work of Moy-Prasad, especially \cite{MP2}, to such nonlinear covers.  Pan claims on occasion that the results of Moy-Prasad generalize ``clearly'' to nonlinear covers; we respectfully disagree and find that a more complete treatment of ``depth'' and ``unrefined minimal $K$-types'' is required for these groups.  The results of this paper should fill in some gaps in \cite{Pan1}.

We mention three purposes for studying the depth zero supercuspidal representations of nonlinear covers of reductive groups.  First, they play a foundational role in the construction of all tame supercuspidal representations, if one believes that the work of Yu \cite{Yu2} extends in some form to the nonlinear case.

Secondly, combined with earlier work of the second author \cite{Wei}, the study of depth zero supercuspidal representations together with a study of metaplectic tori, should lead to a ``local Langlands conjecture'' for metaplectic groups.  Namely, recent work of DeBacker-Reeder \cite{D-R} parameterizes (general position) depth zero representations of reductive groups, in terms of suitable Langlands parameters.  They primarily rely on the work of Moy-Prasad \cite{MP2} (for depth zero representations) and Langlands classification of the representations of (unramified) tori \cite{Lan}.  By generalizing the work of Moy-Prasad to nonlinear covers (the goal of this work), and generalizing Langlands classification to ``metaplectic tori'' (the recent work of \cite{Wei}), the work of DeBacker-Reeder might also be generalized.

The third goal of this work is a more thorough understanding of the supercuspidal representations of the (classical) metaplectic group $Mp_{2n}$.  The first author is particularly interested in the theta correspondence, the study of which requires a thorough knowledge of the representations of the metaplectic group.  In the final section of this paper, we explain this special case in more detail.  On the one hand, our treatment of supercuspidal representations (like that of Joyner \cite{Joy1}) does not provide a local ``Shimura correspondence'' -- it does not directly elucidate connections to representations of orthogonal groups.  On the other hand, it provides a complete classification of genuine, depth zero, supercuspidal irreducible representations of the metaplectic group in odd residue characteristic.  Such a classification should be helpful in future studies of the theta correspondence.

\subsection{Acknowledgements}

The second author would like to thank Stephen DeBacker and Gopal Prasad for helpful discussions during the preparation of this paper.  In addition, he appreciates the hospitality and work-conducive weather at the University of Michigan, during the time when this paper was started.

\section{Structure Theory}

The following notations will be fixed throughout the paper:
\begin{itemize}
\item
$p$ is a prime number, and the field $F$ is a finite extension of $\rats_p$.
\item
$\bar F$ is an algebraic closure of $F$, and $F^{unr}$ is the maximal unramified extension of $F$ in $\bar F$.
\item
$\OO$ denotes the valuation ring of $F$, and $\OO^{unr}$ the valuation ring of $F^{unr}$.
\item
$\ff{f}$ is the residue field of $\OO$, and $\ff{\bar f}$ the residue field of $\OO^{unr}$.
\item
$\Gamma = Gal(F^{unr} / F)$, identified naturally with $Gal(\ff{\bar f} / \ff{f})$.
\end{itemize}

We use the following conventions, throughout the paper:
\begin{itemize}
\item
A boldface letter, such as $\alg{H}$ will be used to denote an algebraic variety over $F$.
\item
An ordinary letter, such as $H$, will be used to denote the $F$-points of an $F$-variety $\alg{H}$.  $H$ will be endowed with the topology arising from the topology on $F$ as a local field.
\item
A letter such as $H^{unr}$ will be used to denote the $F^{unr}$-points of an $F$-variety $\alg{H}$, endowed with the natural action of $\Gamma$.
\item
An underlined boldface letter, such as $\scheme{H}_x$ will be used to denote a scheme over $\OO$.
\item
An ordinary letter (with a subscript), such as $H_x$ will be used to denote the $\OO$-points of an $\OO$-scheme $\scheme{H}_x$, viewed naturally as a topological group.
\item
A letter, such as $H_x^{unr}$ will be used to denote the $\OO^{unr}$-points of an $\OO$-scheme $\scheme{H}_x$, endowed with the natural action of $\Gamma$.
\item
An overlined boldface letter, such as $\alg{\bar H}$ will be used to denote a variety over $\ff{f}$.
\item
An overlined letter, such as $\bar H$, will be used to denote the $\ff{f}$-points of an $\ff{f}$-variety $\alg{\bar H}$.
\end{itemize}

For representations, we use the following conventions:
\begin{itemize}
\item
When $G$ is a topological group, and there exists a basis of neighborhoods of the identity consisting of compact open subgroups, a representation of $G$ will always mean a smooth representation on a complex vector space.
\item
An irreducible smooth representation of $G$ will be called an {\em irrep} of $G$.
\end{itemize}

We fix the following groups and associated data:
\begin{itemize}
\item
The multiplicative group over $F$ will be denoted by $\alg{G}_m$, and over $\ff{f}$ it will be denoted by $\alg{\bar G}_m$.
\item
$\alg{G}$ is a connected reductive algebraic group over the field $F$.  Later, we will place more restrictions on $\alg{G}$.
\item
$\alg{S}$ is a maximal $F$-split torus in $\alg{G}$, and $\alg{T}$ is the centralizer of $\alg{S}$ in $\alg{G}$ (which is a maximal torus in $\alg{G}$, split over $F^{unr}$).
\end{itemize}

\subsection{Groups and buildings}

Our treatment of buildings follows Moy-Prasad \cite{MP2}, and DeBacker \cite{Deb}.  Let $\mc{B}(G)$ denote the (enlarged) Bruhat-Tits building of $G$, identified with the $\Gamma$-fixed points of the building $\mc{B}(G^{unr})$.  If $\alg{S}$ is a maximal $F$-split torus in $\alg{G}$, we write $\mc{A}(S)$ for the corresponding apartment in $\mc{B}(G)$.  Following \cite{Deb}, if $\Omega \subset \mc{A}(S)$, we write $A(\Omega, \mc{A}(S))$ for the smallest affine subspace of $\mc{A}(S)$ containing $\Omega$.

If $x \in \mc{B}(G)$, we write $G_x$ for the parahoric subgroup of $G$ attached to $x$, and $G_x^+$ for its pro-unipotent radical.  We write $G_x^{unr}$ and $G_x^{unr, +}$ for the parahoric subgroup of $G^{unr}$ attached to $x$, and its pro-unipotent radical.

By results of Bruhat-Tits \cite{B-T}, there exists a unique, up to unique isomorphism, smooth group scheme $\scheme{G}_{x}$ over $\OO$, endowed with isomorphisms:
\begin{itemize}
\item
$\scheme{G}_x \times_{\OO} F \isom \alg{G}$.
\item
$\scheme{G}_x(\OO) \isom G_x$.
\end{itemize}
We write $\alg{\bar G}_x$ for the special fibre of $\scheme{G}_x$, and $\alg{\bar M}_x$ for the maximal reductive quotient of $\alg{\bar G}_x$.  There is a canonical identification of groups with $\Gamma$-action:
$$\alg{\bar M}_x(\ff{\bar f}) \isom G_x^{unr} / G_x^{unr, +}.$$
In particular, the $\Gamma$-fixed points yield an isomorphism:
$$\bar M_x = \alg{\bar M}_x(\ff{f}) \isom G_x / G_x^+.$$

\subsection{Deligne-Brylinski Extensions}

Let $F_{Zar}$ denote the big Zariski site, of schemes of finite type over $F$.  Let $\sheaf{K}_{2}$ denote the sheaf of abelian groups on $F_{Zar}$, obtained via Quillen's algebraic K-theory.  We may also view $\alg{G}$ as a sheaf of groups on $F_{Zar}$, by the associated functor of points.  Deligne and Brylinski \cite{D-B} extensively study the category of central extensions of $\alg{G}$ by $\sheaf{K}_2$ in the category of sheaves of groups on $F_{Zar}$.  Let $\alg{G}'$ be such a central extension; in particular, for the unramified extension $F^{unr}/F$, there is then a short exact sequence:
$$1 \rightarrow \sheaf{K}_2(F^{unr}) \rightarrow \alg{G}'(F^{unr}) \rightarrow \alg{G}(F^{unr}) \rightarrow 1.$$
Pushing forward via the tame symbol $\sheaf{K}_2(F^{unr}) \twoheadrightarrow \ff{\bar f}^\times$, one arrives at a central extension of groups:
$$1 \rightarrow \ff{\bar f}^\times \rightarrow \tilde G^{unr} \rightarrow G^{unr} \rightarrow 1.$$
Furthermore, $\Gamma$ acts on the above exact sequence, and the fixed points yield the following central extension (Proposition 12.9 of \cite{D-B}) of topological groups:
$$1 \rightarrow \ff{f}^\times \rightarrow \tilde G \rightarrow G \rightarrow 1,$$
where $\tilde G$ is the pushforward of $\alg{G}'(F)$ via the tame Hilbert symbol $\sheaf{K}_2(F) \twoheadrightarrow \ff{f}^\times$.  We refer to $\tilde G$ as a (tame) {\em nonlinear cover} of the group $G$.

If $H$ is any subgroup of $G$, we write $\tilde H$ for its preimage in $\tilde G$.

\subsection{Residual Extensions}

By results of Section 12.8-12.12 of \cite{D-B}, the following is known:
\begin{prop}
\label{RE}
Suppose that $x \in \mc{B}(G)$.  One may construct a central extension $\alg{\bar M}_x'$ of $\alg{\bar M}_x$ by $\alg{\bar G}_m$ (of reductive groups over $\ff{f}$), from which the central extension in the top row is obtained via pullback from the bottom row, in the following ($\Gamma$-equivariant commutative) diagram:
$$\xymatrix{ 1 \ar[r] & \ff{\bar f}^\times \ar[r] \ar[d]^{=} & \tilde G_x^{unr} \ar[r] \ar[d] & G_x^{unr}  \ar[r] \ar[d] & 1 \\
1 \ar[r] & \ff{\bar f}^\times \ar[r] \ar[d]^{=} & \alg{\bar G}_x'(\ff{\bar f}) \ar[r] \ar[d] & \alg{\bar G}_x(\ff{\bar f}) \ar[r] \ar[d] & 1 \\
1 \ar[r] & \alg{\bar G}_m(\ff{\bar f}) \ar[r] & \alg{\bar M}_x'(\ff{\bar f}) \ar[r] & \alg{\bar M}_x(\ff{\bar f}) \ar[r] & 1.}$$
\end{prop}

In particular, the construction of Deligne-Brylinski yields the following useful corollary:
\begin{cor}
\label{SCSplit}
Suppose that $\alg{\bar M}_x$ is a simply-connected group over $\ff{f}$.  Then, the cover $\tilde G_x^{unr}$ splits ($\Gamma$-equivariantly) over $G_x^{unr}$.
\end{cor}
\proof
If $\alg{\bar M}_x$ is simply-connected, then every extension of $\alg{\bar M}_x$ by $\alg{\bar G}_m$ splits.  Hence the bottom row of the diagram in the previous proposition splits.  Since the top row is obtained via pullback, $\tilde G_x^{unr}$ splits over $G_x^{unr}$.
\qed

The group $\alg{\bar M}_x$ is often, but not always, simply-connected, when $\alg{G}$ is simply-connected.  For example, when $\alg{G} = \alg{G}_2$ (a simply-connected simple split group), the group $\alg{\bar M}_x$ can be isomorphic to $\alg{\bar G}_2$, $\alg{\overline{SL}}_3$, or to $\alg{\overline{SO}}_4$ (Section 10 of \cite{GY1}), for three suitable choices of vertex $x$ in the building.  The first two of these groups are simply-connected, yielding a splitting of the cover.  The third is not simply-connected, and we do not know whether the cover splits in this case.

On the other hand, when $\alg{G} = \alg{Sp}_{2n}$, the groups $\alg{\bar M}_x$ are simply-connected for every vertex $x \in \mc{B}(G)$.  We discuss this case in more detail in the final section of this paper.

\subsection{Splittings}
\label{Spl}
It is important to observe the following elementary fact:
\begin{prop}
Suppose that $H$ is a pro-$p$ subgroup of $G$, for example, a closed subgroup of $G_x^+$ for some point $x$ in the building.  Then, the central extension:
$$1 \rightarrow \ff{f}^\times \rightarrow \tilde H \rightarrow H \rightarrow 1$$
splits uniquely.
\end{prop}
\proof
All topological central extensions of the pro-$p$ group $H$ by the prime-to-$p$ cyclic group $\ff{f}^\times$ split.  Furthermore, the set of splittings is a torsor for the trivial group $Hom(H, \ff{f}^\times)$.
\qed

Suppose that $\alg{U}$ is a unipotent subgroup of $\alg{G}$, defined over $F$.  The central extension of $\alg{G}$ by $\sheaf{K}_2$ then restricts to a central extension:
$$1 \rightarrow \sheaf{K}_2 \rightarrow \alg{U}' \rightarrow \alg{U} \rightarrow 1.$$
From Proposition 3.2 of \cite{D-B}, the above central extension splits uniquely.  Hence, the following exact sequence of topological groups splits canonically and $\Gamma$-equivariantly.:
$$1 \rightarrow \ff{\bar f}^\times \rightarrow \tilde U^{unr} \rightarrow U^{unr} \rightarrow 1.$$
The existence and uniqueness of the topological splitting of the cover $\tilde U \twoheadrightarrow U$ follows from the fact that $U$ is equal to the union if its pro-$p$ subgroups.

\subsection{Conjugation of compact opens}
If $g \in G$, and $x \in \mc{B}(G)$, we write $gx$ for the image of $x$ under the action of $g$ on the building.  We also write $Int[g] \colon G \rightarrow G$ for the action of conjugation by $g$, i.e., inner automorphisms.  Since $\tilde G$ is a central extension of $G$, the action $Int$ of $G$ on $G$ lifts to an action $Int$ of $G$ on $\tilde G$.  Namely, if $g \in G$, and $\tilde h \in \tilde G$, we define $Int[g](\tilde h) = \tilde g \tilde h \tilde g^{-1}$, where $\tilde g$ is any element of $\tilde G$ which lifts $g$.

Given $x \in \mc{B}(G)$, we have a short exact sequence of groups:
$$1 \rightarrow \ff{f}^\times \rightarrow \tilde G_x \rightarrow G_x \rightarrow 1.$$
If $g \in G$, then it is known that: $G_{gx} = Int[g](G_x) = g G_x g^{-1}$.
It follows that the same is true for covers:
$$\tilde G_{gx} = Int[g](\tilde G_x).$$
Similarly, $\tilde G_{gx}^+ = Int[g](\tilde G_x^+)$.

Since there is a unique splitting of the cover $\tilde G_x^+ \rightarrow G_x^+$, it is safe to view $G_x^+$ as a subgroup of $\tilde G$, for every $x \in \mc{B}(G)$.
By the uniqueness of the splitting, it follows that:
\begin{prop}
Suppose that $g \in G$.  Then $Int[g](G_x^+) = G_{gx}^+$ for all $x \in \mc{B}(G)$, where $G_x^+$ and $G_{gx}^+$ are viewed as subgroups of $\tilde G$.
\end{prop}

Using the splitting of $\tilde G_x^+ \rightarrow G_x^+$, and the residual extension (Proposition \ref{RE}), there is a natural isomorphism:
$$\tilde G_x / G_x^+ \isom \bar M_x' = \alg{\bar M}_x'(\ff{f}).$$
\begin{cor}
The isomorphism $Int[g] \colon \tilde G_x \rightarrow \tilde G_{gx}$ descends to an isomorphism:
$$Int[g] \colon \bar M_x' \rightarrow \bar M_{gx}'.$$
\end{cor}
The compatibility with unramified extensions demonstrates that this arises from a $\Gamma$-equivariant isomorphism:
$$Int[g] \colon \alg{\bar M}_x'(\ff{\bar f}) \rightarrow \alg{\bar M}_{gx}'(\ff{\bar f}).$$

\begin{prop}
Suppose that $x,y \in \Build(G)$.  Suppose that the intersection $G_{x} \cap G_{y}$ surjects onto both $\bar M_x$ and $\bar M_{y}$.  Then there is a natural $\Gamma$-equivariant isomorphism:
$$i_{x,y} \colon \alg{\bar M}_x'(\ff{\bar f}) \rightarrow \alg{\bar M}_{y}'(\ff{\bar f}).$$
In particular, $\bar M_x' \isom \bar M_y'$.
\label{IDxy}
\end{prop}
\proof
If $G_x \cap G_{y}$ surjects onto both $\bar M_x$, and $\bar M_{y}$, then Moy-Prasad \cite{MP2} demonstrate that the facets of $\mc{B}(G)$ containing $x$ and $y$ are strongly associated (0-associated, in the sense discussed in Section 3 of \cite{Deb}).  It follows, from Lemma 3.5.1 of \cite{Deb}, that there are equalities:
$$G_x^{unr,+} \cap G_y^{unr} = G_x^{unr,+} \cap G_y^{unr,+} = G_x^{unr} \cap G_y^{unr,+}.$$
Observe that there are natural $\Gamma$-equivariant isomorphisms:
$$\alg{\bar M}_x'(\ff{\bar f}) \isom \frac{\tilde G_x^{unr}}{G_x^{unr,+}} \isom \frac{\tilde G_x^{unr} \cap \tilde G_y^{unr}}{G_x^{unr,+} \cap G_y^{unr}},$$
$$\alg{\bar M}_y'(\ff{\bar f}) \isom \frac{\tilde G_y^{unr}}{G_y^{unr,+}} \isom \frac{\tilde G_y^{unr} \cap \tilde G_x^{unr}}{G_y^{unr,+} \cap G_x^{unr}}.$$
Since $G_x^{unr, +} \cap G_y^{unr} = G_y^{unr, +} \cap G_x^{unr}$, it follows that there is a natural $\Gamma$-equivariant isomorphism:
$$i_{x,y} \colon \alg{\bar M}_x'(\ff{\bar f}) \rightarrow \alg{\bar M}_y'(\ff{\bar f}).$$
\qed

\section{Genuine depth zero representations}
Fix a character $\epsilon \colon \ff{f}^\times \rightarrow \complex^\times$.  Suppose that $H$ is any group, and we are given a central extension of groups:
$$1 \rightarrow \ff{f}^\times \rightarrow H' \rightarrow H \rightarrow 1.$$
Suppose that $(\pi, V)$ is a representation of $H'$ on a complex vector space.  Then, we say that $(\pi, V)$ is $\epsilon$-genuine if, for all $z \in \ff{f}^\times$, and all $v \in V$, $\pi(z) v = \epsilon(z) \cdot v$.

Clearly, $\epsilon$-genuine representations factor through the quotient $H' / Ker(\epsilon)$.  In particular, when $\epsilon$ is a quadratic character ($Im(\epsilon) = \{ \pm 1 \}$), the $\epsilon$-genuine representations factor through a quotient $H'/Ker(\epsilon)$ which is a {\em double-cover} of $H$.  Such representations include the genuine representations of metaplectic groups that appear most frequently in the literature.  These will be discussed in the final section.

\begin{defn}
Suppose that $(\pi, V)$ is an $\epsilon$-genuine irrep of $\tilde G$.  We say that $(\pi, V)$ has {\em depth zero}, if, for some $x \in \mc{B}(G)$, the space $V^{G_x^+}$ is nontrivial.
\end{defn}
It is the goal of this paper to study the depth zero representations of $\tilde G$, following the methods of Moy-Prasad \cite{MP2}.

\subsection{Jacquet modules}

Suppose that $(\pi, V)$ is an $\epsilon$-genuine irrep of $\tilde G$.  Suppose that $\alg{P} = \alg{L} \alg{U}$ is a $F$-parabolic subgroup of $\alg{G}$, with unipotent radical $\alg{U}$ and $F$-Levi subgroup $\alg{L}$.  The canonical splitting (Section \ref{Spl}) of $\tilde U$ over $U$ allows us to identify $U$ (and its opposite, $U^-$) with subgroups of $\tilde G$.  In this way, we may define the Jacquet module:
$$J_U(V) = V / V(U), \mbox{ where } V(U) = span \{ \pi(u) v - v \colon u \in U, v \in V \}.$$
In other words, $J_U(V)$ is the maximal quotient of $V$ on which $U$ acts trivially via $\pi$.

The uniqueness of the splitting of $\tilde U$ over $U$ implies that the splitting is stable under conjugation by $L$.  It follows that the action $\pi$ of $\tilde G$ on $V$ descends to an action $\pi_U$ of $\tilde L$ on $J_U(V)$.  This defines an $\epsilon$-genuine representation of $(\pi_U, J_U(V))$ of $\tilde L$.

Finally, we claim that the result known as Jacquet's Lemma generalizes to nonlinear covers.  Jacquet's Lemma, for linear groups, is stated as Theorem 2.2 in \cite{MP2} and a proof is given by Silberger in \cite{Sil}.  The proof can also be found in the widely circulated notes of Casselman and those of DeBacker.  The proof of Jacquet's Lemma adapts, ``mutatis mutandis'' to nonlinear covers:
\begin{prop}
Suppose that $K$ is an open compact pro-$p$ subgroup of $G$, which has the Iwahori decomposition with respect to $(P,L,U)$:
$$K = (K \cap U^-) \cdot (K \cap L) \cdot (K \cap U).$$
Then, viewing $K$ (in the unique way) as a subgroup of $\tilde G$, the map $V \mapsto J_U(V)$ restricts to a surjective map:
$$V^K \twoheadrightarrow (J_U(V))^{K \cap L}.$$
\end{prop}

\subsection{Depth zero types}
The following definitions are directly adapted from Moy-Prasad \cite{MP2}.
\begin{defn}
An $\epsilon$-genuine {\em depth zero type} for $\tilde G$ is a pair $(x, \sigma)$, where:
\begin{itemize}
\item
$x$ is a point in the building $\mc{B}(G)$.
\item
$\sigma$ is a $\epsilon$-genuine cuspidal irrep of $\bar M_x'$.
\end{itemize}
\end{defn}

\begin{defn}
Suppose that $(x,\sigma)$ and $(y, \tau)$ are two $\epsilon$-genuine depth zero types.  We say that $(x,\sigma)$ and $(y, \tau)$ are {\em associate} if there exists $g \in G$, such that:
\begin{itemize}
\item
$G_x \cap G_{gy}$ surjects onto both $\bar M_x$ and $\bar M_{gy}$,  and hence yields an isomorphism by \ref{IDxy}:
$$i_{x,gy} \colon \bar M_x' \rightarrow \bar M_{gy}'.$$
\item
The representation $\sigma$ of $\bar M_x'$ is isomorphic to the pullback $i_{x,gy}^\ast \tau$ of $\tau$ via the above isomorphism.
\end{itemize}
\end{defn}

\begin{prop}
If $(\pi, V)$ is an $\epsilon$-genuine irrep of $\tilde G$, of depth zero, and if $y \in \mc{B}(G)$ is such that $V^{G_y^+} \neq 0$, then every irreducible $\tilde G_y$-submodule of $V^{G_y^+}$ contains an $\epsilon$-genuine depth zero type of the form $(x, \sigma)$, where $\tilde G_x \subset \tilde G_y$.
\end{prop}
\proof
We adapt Theorem 5.2 of \cite{MP1}, following the proof in Section 7.1 of \cite{MP1}.  Suppose that $\tau$ is an irreducible $\tilde G_y$-submodule of $V^{G_y^+}$.  In particular, $\tau$ is naturally an $\epsilon$-genuine representation of the finite group $\tilde G_y / G_y^+ = \bar M_y'$.  As an irreducible genuine representation of $\bar M_y' = \alg{\bar M}_y'(\ff{f})$, Harish-Chandra's theory applies, and $\sigma$ contains an irreducible cuspidal (necessarily $\epsilon$-genuine) representation $\sigma$ of (the $\ff{f}$-points of) a $\ff{f}$-Levi subgroup $\alg{\bar L}'$ of a parabolic subgroup $\alg{\bar P}'$ of $\alg{\bar M}_y'$.

Every such pair, $\alg{\bar L}' \subset \alg{\bar P}'$, consisting of a Levi subgroup of a parabolic subgroup of $\alg{\bar M}_y'$, arises as the pullback of a pair $(\alg{\bar L}, \alg{\bar P})$, consisting of a Levi subgroup of a parabolic subgroup of $\alg{\bar M}_y$:
$$\xymatrix{
1 \ar[r] & \alg{\bar G}_m \ar[r] & \alg{\bar M}_y' \ar[r] & \alg{\bar M}_y \ar[r] & 1 \\
1 \ar[r] & \alg{\bar G}_m \ar[r] \ar[u] & \alg{\bar P}' \ar[r] \ar[u] & \alg{\bar P} \ar[u] \ar[r] & 1 \\
1 \ar[r] & \alg{\bar G}_m \ar[r] \ar[u] & \alg{\bar L}' \ar[r] \ar[u] & \alg{\bar L} \ar[u] \ar[r] & 1.
}$$

The group $\bar L$ arises as a quotient $G_x / G_x^+$, for some point $x \in \mc{B}(G)$, from which we find that $\bar L'$ arises as a quotient $\tilde G_x / G_x^+$.  Thus, we find that $\tau$ contains the $\epsilon$-genuine depth zero type $(x, \sigma)$.
\qed

\begin{prop}
Suppose that $(x, \sigma)$ and $(y, \tau)$ are two $\epsilon$-genuine depth zero types, occurring in an $\epsilon$-genuine irrep $(\pi, V)$ of $\tilde G$.  Then $(x, \sigma)$ and $(y, \tau)$ are associate.
\end{prop}
\proof
We directly follow the proof in Section 7.2 of Moy-Prasad \cite{MP1} -- their proof holds, with minimal variation.  Note that $\sigma$ is a cuspidal $\epsilon$-genuine irrep of $\tilde G_x / G_x^+$, and $\tau$ is a cuspidal $\epsilon$-genuine irrep of $\tilde G_y / G_y^+$.  Let $V_\sigma$ and $V_\tau$ be subspaces of $V$, on which $\tilde G_x$ and $\tilde G_y$ acts irreducibly via $\sigma$ and $\tau$, respectively.  Let $p_\tau \colon V \rightarrow V_\tau$ be a $\tilde G_y$-equivariant projection onto the subspace $V_\tau$.  By the irreducibility of $\pi$, there exists $\tilde g \in \tilde G$ (lifting $g \in G$), such that the following linear map is nonzero:
$$\phi = p_\tau \circ \pi(g^{-1}) \colon V_\sigma \rightarrow V_\tau.$$
Now, for all $h \in \tilde G_x \cap Int[g](\tilde G_y)$, we find that:
$$\phi \circ \sigma(h) = \tau(g^{-1} h g) \circ \phi.$$

If $G_x \cap G_{gy}$ projects onto $\bar M_x$, and onto $\bar M_{gy}$, then there is an isomorphism:
$$i_{x,gy} \colon \bar M_x' \rightarrow \bar M_{gy}'.$$
Recalling the definition of this isomorphism, $i_{x,gy}$ arises from the equality:
$$\frac{\tilde G_x \cap \tilde G_{gy}}{G_x^+ \cap G_{gy}} = \frac{\tilde G_x \cap \tilde G_{gy}}{G_x \cap G_{gy}^+}.$$
From the equality $\phi \circ \sigma(h) = \tau(g^{-1} h g) \circ \phi$, for all $h \in \tilde G_x \cap \tilde G_{gy}$, we find that $\phi$ determines an isomorphism:
$$i_{x,gy}^\ast \tau \isom \sigma.$$
Hence $(x,\sigma)$ and $(y, \tau)$ are associates, as desired.

If $G_x \cap G_{gy}$ does not project onto $\bar M_x$, then as noted in \cite{MP1}, the intersection $\bar M_x \cap \bar M_{gy}$ projects onto (the $\ff{f}$-points of) a proper parabolic subgroup $\bar P_x$, with ($\ff{f}$-points of) unipotent radical $\bar U_x$.  Furthermore, the preimage of $\bar U_x$ in $\bar M_{gy}'$ contains a nontrivial unipotent subgroup $\bar W \subset \bar M_{gy}'$.  Then, we find that:
\begin{itemize}
\item
$Int[g](\tau)$ contains the trivial representation of $\bar W$.
\item
$\sigma$ agrees with $Int[g](\tau)$ on $\bar W$.
\end{itemize}
But, since $\sigma$ is cuspidal, $\sigma$ cannot be trivial on a unipotent subgroup of $\bar M_{gy}'$, a contradiction.  It follows that $\bar M_x \cap \bar M_{gy}$ projects onto $\bar M_x$.  By the same argument, we find that $\bar M_x \cap \bar M_{gy}$ projects onto $\bar M_{gy}$.
\qed

\begin{prop}
Suppose that $(x, \sigma)$ is an $\epsilon$-genuine depth zero type, occurring in an $\epsilon$-genuine irrep $(\pi, V)$ of $\tilde G$.  If $(y, \tau)$ is an $\epsilon$-genuine depth zero type, and $(y, \tau)$ is associate to $(x, \sigma)$, then $(y, \tau)$ occurs in $(\pi, V)$.
\end{prop}
\proof
We adapt Proposition 6.2 of \cite{MP2} here.  If $(x, \sigma)$ and $(y, \tau)$ are associate, then by replacing $y$ by $g^{-1} y$ for some $g$, we may assume that $G_x \cap G_y$ surjects onto $\bar M_x$ and onto $\bar M_y$.  The structural aspects of Proposition 6.2 of \cite{MP2} carry over; in the language of \cite{Deb}, the facets containing $x$ and $y$ are strongly associated, and we write $B$ for the smallest affine subspace of $\mc{B}(G)$ containing both facets.  Without loss of generality, we may assume that $x$ and $y$ are in adjacent components of $B$.  As in \cite{MP2}, there exists a third point $z \in B$, such that $G_z$ contains both $G_x$ and $G_y$, and such that $\bar M_x$ and $\bar M_y$ are Levi subgroups of associate parabolic subgroups $\bar P_x$ and $\bar P_y$ in $\bar M_z$.  Furthermore, $G_z^+ \subset G_x^+$ and $G_z^+ \subset G_y^+$.  It follows that:
$$\bar M_z' \supset \bar M_x', \mbox{ and } \bar M_z' \supset \bar M_y'.$$

Now, if $(x, \sigma)$ is an $\epsilon$-genuine depth zero type occurring in $(\pi, V)$, then $V^{G_x^+}$ (viewed as a representation of $\bar M_x' = \tilde G_x / G_x^+$) contains the representation $\sigma$ as a constituent.  Since $G_z^+ \subset G_x^+$, we have $V^{G_z^+} \supset V^{G_x^+}$.  Therefore, $V^{G_z^+}$, viewed as a representation of $\bar M_z'$, has an irreducible constituent $\xi$ which arises via parabolic induction from the representation $\sigma$ of $\bar M_x'$.  Since $\bar M_x'$ and $\bar M_y'$ are quotients of associate parabolics, and $(x,\sigma)$ is associate to $(y,\tau)$,  it follows that the restriction of $\xi$ to $\bar M_y'$ contains $\tau$.  Thus, we find that $(y, \tau)$ also occurs in $(\pi, V)$.
\qed

To summarize, we have proven the following:
\begin{thm}
Suppose that $(\pi, V)$ is an $\epsilon$-genuine depth zero representation of $\tilde G$.  Then $(\pi, V)$ contains a depth zero type.  Moreover, if $(x, \sigma)$ and $(y,\tau)$ are two $\epsilon$-genuine depth zero types, and $(x, \sigma)$ occurs in $(\pi, V)$, then $(y, \tau)$ occurs in $(\pi, V)$ if and only if $(y, \tau)$ is associate to $(x, \sigma)$.
\end{thm}

\subsection{Induced representations}
Suppose that $(x, \sigma)$ is an $\epsilon$-genuine depth zero type for $\tilde G$.  Suppose that $G_x$ is a maximal parahoric subgroup of $G$.  Let $\tilde N_x = N_{\tilde G}(\tilde G_x)$ be its normalizer in $\tilde G$.  Let $N_x = N_G(G_x)$ be the normalizer of $G_x$ in $G$, and observe that $\tilde N_x$ is the preimage of $N_x$ in $\tilde G$, as suggested by the notation.  Following Section 6.5 of \cite{MP2}, observe that $N_x$ is compact modulo its center, and so is $\tilde N_x$.  Let $E(\sigma)$ be the set of isomorphism classes of irreps of $\tilde N_x$, which contain $\sigma$ upon restriction to $\tilde G_x$.  We now adapt Proposition 6.6 of \cite{MP2} to the nonlinear case:
\begin{prop}
Suppose that $\tau \in E(\sigma)$.  Let $\pi = c-Ind_{\tilde N_x}^{\tilde G} \tau$ be the compactly-induced representation of $\tilde G$.  Then $\pi$ is an irreducible $\epsilon$-genuine depth zero supercuspidal representation of $\tilde G$.  Moreover, any irreducible admissible representation of $\tilde G$, which contains $\sigma$ upon restriction to $\tilde G_x$, is isomorphic to $c-Ind_{\tilde N_x}^{\tilde G} \tau$, for some $\tau \in E(\sigma)$.
\end{prop}
\proof
The same proof as in \cite{MP2} goes through, almost without changes.  It suffices to prove that the $\tau$-spherical Hecke algebra $H_\tau = H(\tilde G // \tilde N_x, \tau)$ is one-dimensional.  For this, it suffices to prove that any element of this Hecke algebra has support in the ``trivial coset'' $\tilde N_x$.  For this, observe that for $f \in H_\tau$, and $\tilde g \in \tilde G$, a necessary condition for $f(\tilde g) \neq 0$ is that:
$$\tilde G_x \cap \tilde G_{gx} \mbox{ surjects onto } \bar M_x' \mbox{ and onto } \bar M_{gx}'.$$
For this, it is necessary and sufficient to demonstrate that:
$$G_x \cap G_{gx} \mbox{ surjects onto } \bar M_x \mbox{ and } \bar M_{gx}.$$
As shown in the proof of Proposition 6.6 of \cite{MP2}, this implies that $g \in N_x$.  Hence $\tilde g \in \tilde N_x$.  We have proven that $H_\tau$ is one-dimensional, and so $\pi$ is an irrep.  By construction, it is $\epsilon$-genuine.  Since the matrix coefficients of $\pi$ are compactly-supported, $\pi$ is supercuspidal.  By Frobenius reciprocity, $\pi$ contains the type $(x, \sigma)$.

The last portion of the proof of Proposition 6.6 of \cite{MP2} goes through without changes, and we find that any irrep $\tau$ of $\tilde G$, which contains $\sigma$ upon restriction to $\tilde G_x$, is isomorphic to $c-Ind_{\tilde N_x}^{\tilde G} \tau$, for some $\tau \in E(\sigma)$.
\qed

One may twist $\epsilon$-genuine irreps $(\pi, V)$ of $\tilde G$ by unramified quasicharacters of $G$; namely, if $\chi \colon G \rightarrow \complex^\times$ is an unramified quasicharacter, then it pulls back to a quasicharacter of $\tilde G$.  Then, $\pi \otimes \chi$ is still an $\epsilon$-genuine irrep of $\tilde G$.  Similarly, we may use $\chi$ to twist an $\epsilon$-genuine representation $\tau$ of $\tilde N_x$.

As in \cite{MP2}, we find that:
$$\left( c-Ind_{\tilde N_x}^{\tilde G} \tau \right) \otimes \chi \isom c-Ind_{\tilde N_x}^{\tilde G} (\tau \otimes \chi).$$

\begin{thm}
Every irreducible, $\epsilon$-genuine, depth zero, supercuspidal representation of $\tilde G$ arises as $c-Ind_{\tilde N_x}^{\tilde G} \tau$, for some depth zero, $\epsilon$-genuine type $(x,\sigma)$, and some $\tau \in E(\sigma)$.
\end{thm}
\proof
The proof of Propositions 6.7 in \cite{MP2} goes through without changes, since it only relies on subgroups of $G$ which split uniquely in $\tilde G$.  The proof of Proposition 6.8 of \cite{MP2} then proceeds without changes.
\qed

\section{The Metaplectic Group}

In this section, we consider the group $\alg{G} = \alg{Sp}_{2n}$, and its universal central extension $\alg{G}'$ by $\sheaf{K}_2$.  This is called Matsumoto's central extension, in (0.2) of \cite{D-B}.  The universal central extension yields a central extension of topological groups:
$$1 \rightarrow \ff{f}^\times \rightarrow \tilde G \rightarrow G \rightarrow 1.$$
Let $\varpi$ denote a uniformizing element of $F$.  We also assume that $p \neq 2$ ($F$ has odd residual characteristic).  Thus, there exists a unique surjective homomorphism $\epsilon$ from $\ff{f}^\times$ to the group $\mu_2 = \{ \pm 1 \}$ of order two.  Pushing forward via this homomorphism yields the classical {\em metaplectic group}:
$$1 \rightarrow \mu_2 \rightarrow Mp_{2n} \rightarrow G \rightarrow 1.$$

It follows that the $\epsilon$-genuine representations of $\tilde G$ correspond precisely to the {\em genuine} representations of $Mp_{2n}$ (those for which $\mu_2$ acts via the faithful sign character).

\begin{prop}
Suppose that $x$ is a vertex in the building of $G = Sp_{2n}$.  Then $\tilde G_x$ splits over $G_x$, and thus the cover $Mp_{2n}$ splits over the maximal parahoric subgroup $G_x$ of $Sp_{2n}$.
\end{prop}
\proof
Let $V$ denote a $2n$-dimensional $F$ vector space, endowed with a nondegenerate symplectic form $\langle \cdot, \cdot \rangle$, such that $G = Sp(V)$.  As discussed by Yu \cite{Yu3}, there is a natural correspondence between the vertices in $\mc{B}(G)$ and the set of $\OO$-lattices $\Lambda$ in $V$, such that:
$$\varpi \Lambda^\perp \subset \Lambda \subset \Lambda^\perp.$$
For such a lattice $\Lambda$, the following are naturally defined $\ff{f}$-vector spaces:
$$W_1 = \Lambda / \varpi \Lambda^\perp, \mbox{ and } W_2 = \Lambda^\perp / \Lambda.$$
The $F$-valued symplectic form on $V$ descends to a natural nondegenerate $\ff{f}$-valued symplectic form on each of these $\ff{f}$-vector spaces.  In \cite{Yu3}, Yu observes that the reductive group $\alg{\bar M}_x$ at $x$ can be identified:
$$\alg{\bar M}_x \isom \alg{Sp}(W_1) \times \alg{Sp}(W_2).$$
One can deduce this from the characterization of $\mc{B}(G)$ within $\mc{B}(SL_{2n})$, as discussed by Kim-Moy \cite{K-M}.  In particular, $\alg{\bar M}_x$ is a simply-connected algebraic group over $\ff{f}$.  The residual extension:
$$1 \rightarrow \alg{\bar G}_m \rightarrow \alg{\bar M}_x' \rightarrow \alg{\bar M}_x \rightarrow 1$$
splits uniquely.  This unique splitting (following Corollary \ref{SCSplit}) splits the cover $\tilde G_x$ of $G_x$.  It directly follows that $Mp_{2n}$ splits over the maximal parahoric $G_x$.
\qed

The normalizer of a maximal parahoric $G_x$ in $G = Sp_{2n}$ is simply $G_x$, since $\alg{G}$ is simply-connected. It follows that:
\begin{prop}
If $x$ is a vertex in $\mc{B}(G)$, then every $\epsilon$-genuine, depth zero type for $\tilde G$ has the form $(x, \epsilon \boxtimes \sigma)$, where $\sigma$ is a cuspidal irrep of $\bar M_x$, and $\bar M_x'$ is identified with $\ff{f}^\times \times \bar M_x'$.
\end{prop}

Let $sgn$ denote the unique nontrivial character of $\mu_2$.  From the previous section, we find the following:
\begin{thm}
Suppose that $(\pi, V)$ is a genuine, irreducible, depth zero supercuspidal representation of $Mp_{2n}$.  Then, there exists a vertex $x$ in $\mc{B}(Sp_{2n})$, and a cuspidal irrep $\sigma$ of $\bar M_x$ (inflated to $G_x$), such that:
$$\pi \isom c-Ind_{\tilde G_x}^{\tilde G} (sgn \boxtimes \sigma).$$

Conversely, given any vertex $x \in \mc{B}(Sp_{2n})$, and any cuspidal irrep $\sigma$ of $\bar M_x$ (inflated to $G_x$), define the representation:
$$\pi = c-Ind_{\tilde G_x}^{\tilde G} (sgn \boxtimes \sigma).$$
Then $\pi$ is a genuine, irreducible, depth zero supercuspidal representation of $Mp_{2n}$.
\end{thm}

Recalling that $\bar M_x$ is a product of symplectic groups over a finite field, this yields a complete constructive description of the genuine depth zero supercuspidal irreps of the metaplectic group, in odd residue characteristic.  Also, observe that these genuine depth zero supercuspidal irreps of the metaplectic group are in natural bijection with the depth zero supercuspidal irreps of the symplectic group, under this assumption of odd residue characteristic.  This generalizes results of Joyner \cite{Joy1}, \cite{Joy2}, from $SL_2$ to $Sp_{2n}$.

\bibliographystyle{amsplain}
\bibliography{MetaplecticDepth0}
\end{document}